\documentclass[psamsfonts,reqno,11pt]{amsart}
\usepackage{amssymb}
\usepackage[all]{xy}

\theoremstyle{plain}

\theoremstyle{definition}

\theoremstyle{remark}

\DeclareSymbolFont{AMSb}{U}{msb}{m}{n}

\DeclareMathSymbol{\A}{\mathbin}{AMSb}{"41}
\DeclareMathSymbol{\B}{\mathbin}{AMSb}{"42}
\DeclareMathSymbol{\C}{\mathbin}{AMSb}{"43}
\DeclareMathSymbol{\D}{\mathbin}{AMSb}{"44}
\DeclareMathSymbol{\E}{\mathbin}{AMSb}{"45}
\DeclareMathSymbol{\F}{\mathbin}{AMSb}{"46}
\DeclareMathSymbol{\G}{\mathbin}{AMSb}{"47}
\DeclareMathSymbol{\HH}{\mathbin}{AMSb}{"48}
\DeclareMathSymbol{\I}{\mathbin}{AMSb}{"49}
\DeclareMathSymbol{\N}{\mathbin}{AMSb}{"4E}
\DeclareMathSymbol{\PP}{\mathbin}{AMSb}{"50}
\DeclareMathSymbol{\Q}{\mathbin}{AMSb}{"51}
\DeclareMathSymbol{\R}{\mathbin}{AMSb}{"52}
\DeclareMathSymbol{\SSS}{\mathbin}{AMSb}{"53}
\DeclareMathSymbol{\T}{\mathbin}{AMSb}{"54}
\DeclareMathSymbol{\U}{\mathbin}{AMSb}{"55}
\DeclareMathSymbol{\V}{\mathbin}{AMSb}{"56}
\DeclareMathSymbol{\W}{\mathbin}{AMSb}{"57}
\DeclareMathSymbol{\X}{\mathbin}{AMSb}{"58}
\DeclareMathSymbol{\Y}{\mathbin}{AMSb}{"59}
\DeclareMathSymbol{\Z}{\mathbin}{AMSb}{"5A}

\def\th#1{\noindent{\bf #1}\bgroup\it}
\def\endth{\egroup\par}

\newcommand{\lleq}{\leqslant}
\newcommand{\ggeq}{\geqslant}

\newcommand{\si}{\sum^\infty_{i=1}}

\newcommand{\sni}{\sum^n_{i=1}}

\newcommand{\smk}{\sum^m_{k=1}}
\newcommand{\smi}{\sum^m_{i=1}}

\newcommand{\LL}{\mathcal L}
\newcommand{\Lr}{{\mathcal L}^r}

\newcommand{\p}{\mathcal P}

\newcommand{\pr}{{\mathcal P}^r}

\newcommand{\ex}{\otimes_nX}

\newcommand{\esx}{\otimes_{n,s}X}

\newcommand{\e}{\otimes_nE}
\newcommand{\ep}{\hat{\otimes}_{n,|\pi|}E}
\newcommand{\es}{\otimes_{n,s}E}
\newcommand{\esp}{\hat{\otimes}_{n,s,|\pi|}E}

\newcommand{\epe}{\hat{\otimes}_{n,\pi}E}
\newcommand{\espe}{\hat{\otimes}_{n,s,\pi}E}

\newcommand{\dpe}{\Delta(\hat{\otimes}_{n,\pi}E)}
\newcommand{\dsp}{\Delta(\hat{\otimes}_{n,s,\pi}E)}
\newcommand{\dpp}{\Delta(\hat{\otimes}_{n,|\pi|}E)}
\newcommand{\dspp}{\Delta(\hat{\otimes}_{n,s,|\pi|}E)}

\newcommand{\iepx}{\check{\otimes}_{n,\epsilon}X}
\newcommand{\iespx}{\check{\otimes}_{n,s,\epsilon}X}
\newcommand{\iep}{\check{\otimes}_{n,|\epsilon|}E}
\newcommand{\iesp}{\check{\otimes}_{n,s,|\epsilon|}E}
\newcommand{\iepe}{\check{\otimes}_{n,\epsilon}E}
\newcommand{\iespe}{\check{\otimes}_{n,s,\epsilon}E}

\newcommand{\idpe}{\Delta(\check{\otimes}_{n,\epsilon}E)}
\newcommand{\idsp}{\Delta(\check{\otimes}_{n,s,\epsilon}E)}
\newcommand{\idpp}{\Delta(\check{\otimes}_{n,|\epsilon|}E)}
\newcommand{\idspp}{\Delta(\check{\otimes}_{n,s,|\epsilon|}E)}
\newcommand{\idpx}{\Delta(\check{\otimes}_{n,\epsilon}X)}
\newcommand{\idspx}{\Delta(\check{\otimes}_{n,s,\epsilon}X)}

\begin{document}

\title[]
{Diagonals of Injective Tensor Products of Banach Lattices with Bases}

\author{Donghai Ji}
\address{Department of Mathematics, Harbin University of Science and Technology, Harbin 150080, China}
\email{jidonghai@126.com}

\author{Byunghoon Lee}
\address{Department of Mathematics, University of Mississippi, University, MS 38677, USA}
\email{blee4@olemiss.edu}

\author{Qingying Bu}
\address{Department of Mathematics, University of Mississippi, University, MS 38677, USA}
\email{qbu@olemiss.edu}

\thanks{The first author is supported by the Heilongjiang NSF (No. A201011) and China NNSF (No. 11171082).}

\date{\today}

\begin{abstract}
Let $E$ be a Banach lattice with a 1-unconditional basis $\{e_i: i \in \N\}$. Denote by $\idpe$ (resp. $\idsp$) the
main diagonal space of the $n$-fold full (resp. symmetric) injective Banach space tensor product, and denote by $\idpp$ (resp. $\idspp$) the main diagonal space of the $n$-fold full (resp. symmetric) injective Banach lattice tensor product.  We show that these four main diagonal spaces are pairwise isometrically isomorphic. We also show that the tensor diagonal $\{e_i\otimes\cdots\otimes e_i: i \in \N\}$ is a 1-unconditional basic sequence in both $\iepe$ and $\iespe$.  
\end{abstract}

\maketitle

\vspace{.2in}
\noindent
{\bf 1. Introduction}
\vspace{.2in}

Let $X$ be a Banach space with a 1-unconditional basis $\{e_i: i \in \N\}$. Then $\{e_{i_1}\otimes\cdots\otimes e_{i_n}:\, (i_1, \dots, i_n) \in \N^n\}$ is a basis of both the $n$-fold full injective tensor product $\iepx$ and the $n$-fold full projective tensor product $\hat{\otimes}_{n,\pi}X$ (see, e.g., \cite{Ge, GR1}), and $\{e_{i_1}\otimes_s\cdots\otimes_s e_{i_n}:\, (i_1, \dots, i_n) \in \N^n, i_1 \ggeq\cdots\ggeq i_n\}$ is a basis of both the $n$-fold symmetric injective tensor product $\iespx$ and the $n$-fold symmetric projective tensor product $\hat{\otimes}_{n,s,\pi}X$ (see, e.g., \cite{GR1}). However, they are not necessary unconditional bases (see, e.g., \cite{KP, Pi, Sch, DDGM, PV, CG}). In particular, the tensor diagonal $\{e_i\otimes\cdots\otimes e_i: i \in \N\}$ is a 1-unconditional basic sequence in both $\hat{\otimes}_{n,\pi}X$ (see, e.g., \cite{Ho, Ca}) and $\hat{\otimes}_{n,s,\pi}X$ (see, e.g., \cite{BB1}); and Holub \cite{Ho} showed that the tensor diagonal $\{e_i\otimes e_i: i \in \N\}$ is a 1-unconditional basic sequence in $\check{\otimes}_{2,\epsilon}X$. By using Holub's method, it is easy to show that the tensor diagonal $\{e_i\otimes\cdots\otimes e_i: i \in \N\}$ is a 1-unconditional basic sequence in $\check{\otimes}_{n,\epsilon}X$. However, Holub's method does not work for $\iespx$. In this paper, by using a result obtained in \cite{BB1} we show that the tensor diagonal $\{e_i\otimes\cdots\otimes e_i: i \in \N\}$ is a 1-unconditional basic sequence in $\iespx$ and we also show that the diagonal projections on both $\iepx$ and $\iespx$ are contractive.

From the positivity perspective, let $E$ be a Banach lattice with a 1-unconditional basis $\{e_i: i \in \N\}$. Then $\{e_{i_1}\otimes\cdots\otimes e_{i_n}:\, (i_1, \dots, i_n) \in \N^n\}$ is a 1-unconditional basis of  the $n$-fold full positive projective tensor product $\hat{\otimes}_{n,|\pi|}E$, and $\{e_{i_1}\otimes_s\cdots\otimes_s e_{i_n}:\, (i_1, \dots, i_n) \in \N^n, i_1 \ggeq\cdots\ggeq i_n\}$ is a 1-unconditional basis of the $n$-fold positive symmetric projective tensor product $\hat{\otimes}_{n,s,|\pi|}E$ (see, e.g., \cite{can, BB1}). In this paper we show that $\{e_{i_1}\otimes\cdots\otimes e_{i_n}:\, (i_1, \dots, i_n) \in \N^n\}$ is a 1-unconditional basis of  the $n$-fold full positive injective tensor product $\iep$, and $\{e_{i_1}\otimes_s\cdots\otimes_s e_{i_n}:\, (i_1, \dots, i_n) \in \N^n, i_1 \ggeq\cdots\ggeq i_n\}$ is a 1-unconditional basis of the $n$-fold positive symmetric injective tensor product $\iesp$.

Let $\dpe$ (resp. $\dsp$) denote the closed subspace generated by the tensor diagonals $\{e_i\otimes\cdots\otimes e_i: i \in \N\}$ in $\epe$ (resp. $\espe$), and let $\dpp$ (resp. $\dspp$) denote the closed sublattice generated by the tensor diagonal $\{e_i\otimes\cdots\otimes e_i: i \in \N\}$ in $\ep$ (resp. $\esp$). Bu and Buskes \cite{BB1} showed that these four main diagonal spaces $\dpe, \dsp, \dpp$, and $\dspp$ are isometrically isomorphic. Now let $\idpe$ (resp. $\idsp$) denote the closed subspace generated by the tensor diagonal $\{e_i\otimes\cdots\otimes e_i: i \in \N\}$ in $\iepe$ (resp. $\iespe$), and let $\idpp$ (resp. $\idspp$) denote the closed sublattice generated by the tensor diagonals $\{e_i\otimes\cdots\otimes e_i: i \in \N\}$ in $\iep$ (resp. $\iesp$). In this paper we use the 1-unconditionality of the tensor diagonal $\{e_i\otimes\cdots\otimes e_i: i \in \N\}$ to show that these four main diagonal spaces $\idpe, \idsp, \idpp$, and $\idspp$ are isometrically isomorphic.

\vspace{.3in}
\noindent
{\bf 2. Preliminaries}
\vspace{.2in}

For a Banach space $X$, let $X^\ast$ denote its topological dual and $B_X$ denote its closed unit ball. For Banach spaces $X_1, \dots, X_n$, and $X$, $Y$, let $\LL(X_1, \dots, X_n; Y)$ denote the space of all continuous $n$-linear operators from $X_1\times \cdots \times X_n$ to $Y$, and $\p(^nX; Y)$ denote the space of all continuous $n$-homogeneous polynomials from $X$ to $Y$.
Let $X_1\otimes\cdots\otimes X_n$ denote the {\it n-fold algebraic tensor product} of $X_1, \dots, X_n$. The {\it injective tensor norm} on $X_1\otimes\cdots\otimes X_n$ is defined by
$$
\|u\|_\epsilon = \sup\left\{\left|\sum_{k=1}^m x^\ast_1(x_{1,k}) \cdots x^\ast_n(x_{n,k})\right|: u = \sum_{k=1}^m x_{1,k}\otimes \cdots\otimes x_{n,k}, x^\ast_i \in B_{X^\ast_i}, 1 \lleq i \lleq n \right\}
$$
for every $u \in X_1\otimes\cdots\otimes X_n$. The completion of $X_1\otimes\cdots\otimes X_n$ with respect to this norm is denoted by $X_1\check{\otimes}_\epsilon \cdots \check{\otimes}_\epsilon X_n$ and called the {\it n-fold injective tensor product} of $X_1, \dots, X_n$. For each $u \in X_1\otimes\cdots\otimes X_n$, say, $u = \sum_{k=1}^m x_{1,k}\otimes \cdots\otimes x_{n,k}$, define $T_u:X_1^\ast\times\cdots\times X_n^\ast \rightarrow \R$ by
$$
T_u(x_1^\ast, \dots, x_n^\ast) = \smk x^\ast_1(x_{1,k}) \cdots x^\ast_n(x_{n,k}), \quad \forall\; x^\ast_i \in X^\ast_i,\; i = 1, \dots, n.
\eqno{(2.1)}
$$
Then $T_u$ is a $n$-linear operator (which does not depend on the representations of $u$) and $T_u \in \LL(X^\ast_1, \dots, X_n^\ast; \R)$ with $\|T_u\| = \|u\|_\epsilon$.

\vspace{.1in}
The following lemma is straightforward from the definition.

\vspace{.1in}
\th{Lemma 2.1.}
Let $T_i: X_i \to X_i$ be bounded linear operators for $i = 1, \dots, n$. Then
$$
\left\|\sum_{k=1}^m T_1(x_{1,k})\otimes \cdots\otimes T_n(x_{n,k})\right\|_\epsilon \lleq \|T_1\| \cdots \|T_n\| \cdot \left\|\sum_{k=1}^m x_{1,k}\otimes \cdots\otimes x_{n,k}\right\|_\epsilon
$$
for every $x_{1,k} \in X_1, \dots, x_{n,k} \in X_n$, $k = 1, \dots, m$.
\endth
\vspace{.1in}

If $X_1 = \cdots = X_n = X$, we write $X_1\otimes\cdots\otimes X_n$ by $\ex$ and  $X_1\hat{\otimes}_\epsilon \cdots \hat{\otimes}_\epsilon X_n$ by $\iepx$. For $x_1\otimes\cdots\otimes x_n \in \ex$, let $x_1\otimes_s\cdots\otimes_s x_n$ denote its symmetrization, that is,
$$
x_1\otimes_s\cdots\otimes_s x_n = \frac{1}{n!}\sum_{\sigma \in \pi(n)} x_{\sigma(1)}\otimes\cdots\otimes x_{\sigma(n)},
\eqno{(2.2)}
$$
where $\pi(n)$ is the group of permutations of $\{1, \dots, n\}$. Then (see, e.g., \cite{Fl}) 
$$
x_1\otimes_s\cdots\otimes_s x_n = \frac{1}{2^nn!}\sum_{\delta_i=\pm 1} \delta_1\cdots \delta_n \Big(\sni \delta_ix_i\Big)\otimes\cdots\otimes \Big(\sni \delta_ix_i\Big).
\eqno{(2.3)}
$$
Let $\esx$ denote the {\it n-fold symmetric algebraic tensor product} of $X$, that is, the linear span of $\{x_1\otimes_s\cdots\otimes_s x_n: x_1, \dots, x_n \in X\}$ in $\ex$. Each $u \in \esx$ has a representation $u = \smk \lambda_k x_{k}\otimes \cdots\otimes x_{k}$ where $\lambda_1, \dots, \lambda_m$ are scalars and $x_1, \dots, x_m$ are vectors in $X$. By linearly extending, (2.2) defines a linear projection $s: \ex \rightarrow \esx$.

\vspace{.1in}
The {\it symmetric injective tensor norm} on $\esx$ is defined by
$$
\|u\|_{s,\epsilon} = \sup\left\{\left|\smk \lambda_k\cdot \Big(x^\ast(x_k)\Big)^n\right|: u = \smk \lambda_k x_{k}\otimes \cdots\otimes x_{k}, x^\ast \in B_{X^\ast} \right\}
$$
for every $u \in \esx$. The completion of $\esx$ with respect to this norm is denoted by $\iespx$ and called the {\it $n$-fold symmetric injective tensor product} of $X$. The linear projection $s: \ex \rightarrow \esx$ can be extended to $\iepx$ with $\iespx$ as its range and for every $u \in \iepx$ (see, e.g., \cite{Fl}), we have
$$
\big\|s(u)\big\|_{s,\epsilon} \lleq \big\|s(u)\big\|_{\epsilon} \lleq \frac{n^n}{n!}\big\|s(u)\big\|_{s,\epsilon}.
\eqno{(2.4)}
$$

\vspace{.1in}
For each $u \in \esx$, say, $u = \smk \lambda_k x_{k}\otimes \cdots\otimes x_{k}$, define $P_u:X^\ast \rightarrow \R$ by
$$
P_u(x^\ast) = \smk \lambda_k\cdot \Big(x^\ast(x_k)\Big)^n, \quad \forall\; x^\ast \in X^\ast.
\eqno{(2.5)}
$$
Then $P_u$ is a $n$-homogeneous polynomial (which does not depend on the representations of $u$) and  $P_u \in \p(^nX^\ast; \R)$ with $\|P_u\| = \|u\|_{s,\epsilon}$.

\vspace{.1in}
The following lemma is straightforward from the definition.

\vspace{.1in}
\th{Lemma 2.2.}
Let $T: X \to X$ be a bounded linear operator. Then
$$
\left\|\smk \lambda_k T(x_{k})\otimes \cdots\otimes T(x_{k})\right\|_{s,\epsilon} \lleq \|T\|^n \cdot \left\|\smk \lambda_k x_{k}\otimes \cdots\otimes x_{k}\right\|_{s,\epsilon}
$$
for every $x_{k} \in X$, $k = 1, \dots, m$.
\endth
\vspace{.1in}

For the basic knowledge about the (symmetric) injective tensor products $\iepx$ and $\iespx$, we refer to \cite{DF, Di, Fl, Mu, R}.

\vspace{.2in}
For Banach lattices $E, E_1, \dots, E_n, F$ with $F$ Dedekind complete, let $\Lr(E_1, \dots, E_n; F)$ denote the Banach lattice of all regular $n$-linear operators from $E_1\times\cdots\times E_n$ to $F$ with its regular operator norm $\|T\|_r = \|\,|T|\,\|$, and let $\pr(^nE;F)$ denote the Banach lattice of all regular $n$-homogeneous polynomials from $E$ to $F$ with its regular norm $\|P\|_r = \|\,|P|\,\|$ (see, e.g., \cite{BB}).

\vspace{.1in}
\th{Lemma 2.3.}
Let $E_1, \dots, E_n, F$ be Banach lattices with $F$ Dedekind complete. Then for any $T \in \Lr(E_1, \dots, E_n;F)$,
\begin{align*}
\|T\|_r  = \inf\Big\{\|S\|:\; & S \in \Lr(E_1, \dots, E_n;F)^+, \\
& |T(x_1, \dots, x_n)| \lleq S(|x_1|, \dots, |x_n|), \; \forall\; x_1 \in E_1, \dots, \;\forall\; x_n \in E_n \Big\}.
\tag{2.6}
\end{align*}
Moreover, $\|T\| \lleq \|T\|_r$.
\endth

\begin{proof}
Let 
\begin{align*}
\mathcal{A} := \Big\{S:\; & S \in \Lr(E_1, \dots, E_n;F)^+,\\
 &  |T(x_1, \dots, x_n)| \lleq S(|x_1|, \dots, |x_n|), \; \forall\; x_1 \in E_1, \dots, \;\forall\; x_n \in E_n \Big\},
\end{align*}
and let $a = \inf\{\|S\|: S \in \mathcal{A}\}$. For any $S \in \mathcal{A}$, it follows from (2.10) in \cite{BB} that $|T| \lleq S$. Thus $\|\,|T|\,\| \lleq \|S\|$ and hence, $\|T\|_r \lleq a$. On the other hand, $|T(x_1, \dots, x_n)| \lleq |T|(|x_1|, \dots, |x_n|)$ for every $x_1 \in E_1, \dots, x_n \in E_n$. Thus $|T| \in \mathcal{A}$ and hence, $a \lleq \|\,|T|\,\| = \|T\|_r$. Therefore, $\|T\|_r = a$ and then it is easy to see that $\|T\| \lleq \|T\|_r$.
\end{proof}

\vspace{.1in}
Similarly, we have the following lemma.

\vspace{.1in}
\th{Lemma 2.4.}
Let $E$ and $F$ be Banach lattices with $F$ Dedekind complete. Then for any $P \in \pr(^nE;F)$,
$$
\|P\|_r  = \inf\Big\{\|R\|:\;  R \in \pr(^nE;F)^+, 
 |P(x)| \lleq S(|x|), \; \forall\; x \in E \Big\}.
\eqno{(2.7)}
$$
Moreover, $\|P\| \lleq \|P\|_r$.
\endth

\vspace{.1in}
For any $u \in E_1\otimes\cdots\otimes E_n$, it is easy to see that the operator $T_u:E_1^\ast\times\cdots\times E_n^\ast \rightarrow \R$ defined in (2.1) is a regular $n$-linear operator and hence, $T_u \in \Lr(E^\ast_1, \dots, E_n^\ast; \R)$. Let $E_1\check{\otimes}_{|\epsilon|}\cdots\check{\otimes}_{|\epsilon|}E_n$ denote the closed sublattice generated by $E_1\otimes\cdots\otimes E_n$ in $\Lr(E^\ast_1, \dots, E_n^\ast; \R)$, called the {\it $n$-fold positive injective tensor product} of $E_1, \dots, E_n$. The norm on $E_1\check{\otimes}_{|\epsilon|}\cdots\check{\otimes}_{|\epsilon|}E_n$ is denoted by $\|\cdot\|_{|\epsilon|}$, that is, for every $u \in E_1\otimes\cdots\otimes E_n$, $\|u\|_{|\epsilon|} = \|T_u\|_r$. By Lemma 2.3, $\|u\|_\epsilon \lleq \|u\|_{|\epsilon|}$. In particular, if $u$ is a positive element in $E_1\otimes\cdots\otimes E_n$, then
$$
\|u\|_{|\epsilon|} = \sup\left\{\left|\sum_{k=1}^m x^\ast_1(x_{1,k}) \cdots x^\ast_n(x_{n,k})\right|: u = \sum_{k=1}^m x_{1,k}\otimes \cdots\otimes x_{n,k}, x^\ast_i \in B^+_{E^\ast_i}, 1 \lleq i \lleq n \right\}.
$$

If $E_1 = \cdots = E_n = E$, we write $E_1\check{\otimes}_{|\epsilon|}\cdots\check{\otimes}_{|\epsilon|}E_n$ by $\iep$.
For any $u \in \es$, it is easy to see that the polynomial $P_u:E^\ast \rightarrow \R$ defined in (2.5) is a regular $n$-homogeneous polynomial and hence, $P_u \in \pr(^nE^\ast; \R)$. Let $\iesp$ denote the closed sublattice generated by $\es$ in $\pr(^nE^\ast; \R)$, called the {\it $n$-fold positive symmetric injective tensor product} of $E$. The norm on $\iesp$ is denoted by $\|\cdot\|_{s,|\epsilon|}$, that is, for every $u \in \es$, $\|u\|_{s,|\epsilon|} = \|P_u\|_r$. By Lemma 2.4, $\|u\|_{s,\epsilon} \lleq \|u\|_{s,|\epsilon|}$. In particular, if $u$ is a positive element in $\es$, then
$$
\|u\|_{s,|\epsilon|} = \sup\left\{\left|\smk \lambda_k\cdot \Big(x^\ast(x_k)\Big)^n\right|: u = \smk \lambda_k x_{k}\otimes \cdots\otimes x_{k}, x^\ast \in B^+_{E^\ast} \right\}.
$$

\vspace{.3in}
\noindent
{\bf 3. Diagonals of Injective Tensor Products}
\vspace{.2in}

In this section we assume that $X$ is a Banach space with a 1-unconditional basis $\{e_i:\, i \in \N\}$. Gelbaum and Lamadrid \cite{Ge} showed that $\{e_i\otimes e_j:\, (i,j) \in \N^2\}$ with the square order is a  basis of $\check{\otimes}_{2,\epsilon}X$. In general, Grecu and Ryan \cite{GR1} showed that $\{e_{i_1}\otimes\cdots\otimes e_{i_n}:\, (i_1, \dots, i_n) \in \N^n\}$ with the order defined in \cite{GR1} is a basis of $\iepx$. They also showed that $\{e_{i_1}\otimes_s\cdots\otimes_s e_{i_n}:\, (i_1, \dots, i_n) \in \N^n, i_1 \ggeq \cdots \ggeq i_n\}$ with the order defined in \cite{GR1} is a basis of $\iespx$.

Let $\idpx$ (resp. $\idspx$) denote the {\it main diagonal space} of $\iepx$ (resp. $\iespx$), that is, the closed subspace spanned in $\iepx$ (resp. in $\iespx$) by the {\it tensor diagonal}  $\{e_i\otimes\cdots\otimes e_i: \, i \in \N\}$. It is known that $\{e_{i_1}\otimes\cdots\otimes e_{i_n}:\, (i_1, \dots, i_n) \in \N^n\}$ and $\{e_{i_1}\otimes_s\cdots\otimes_s e_{i_n}:\, (i_1, \dots, i_n) \in \N^n, i_1 \ggeq \cdots \ggeq i_n\}$, respectively, is not necessary an unconditional basis of $\iepx$ and $\iespx$ (see, e.g., \cite{KP, Pi, Sch, DDGM, PV, CG}). Next we will use the  following Rademacher averaging formula (see, e.g., \cite[Lemma 2.22]{R}) to show that the tensor diagonal $\{e_i\otimes\cdots\otimes e_i:\, i \in \N\}$ is an unconditional basis of both $\idpx$ and $\idspx$, and their diagonal projections are contractive.

\vspace{.1in}
\th{Rademacher Averaging.}
Let $Z_1, \dots, Z_n$ be vector spaces and $x_{i,k} \in Z_i$ for $i = 1, \dots, n$ and $k = 1, \dots, m$. Then
$$
\smk x_{1,k}\otimes\cdots\otimes x_{n,k} = \int_0^1 \left(\smk r_k(t)x_{1,k}\right)\otimes\cdots\otimes \left(\smk r_k(t)x_{n,k}\right)\,dt,
$$
where $\{r_k(t)\}_1^\infty$ is the sequence of Rademacher functions on $[0,1]$.
\endth

\vspace{.1in}
\th{Lemma 3.1.}
The tensor diagonal $\{e_i\otimes\cdots\otimes e_i:\, i \in \N\}$ is a $1$-unconditional basis of $\idpx$ and the projection $Q: \iepx\rightarrow \idpx$ defined by
$$
Q(e_{i_1}\otimes\cdots\otimes e_{i_n}) = \left\{
\begin{array}{ll}
e_{i_1}\otimes\cdots\otimes e_{i_n}, & \mbox{if}\;\; i_1 = \cdots = i_n,\\
0, & \mbox{otherwise}.
\end{array}
\right.
$$
is bounded  with $\|Q\| \lleq 1$.
\endth

\begin{proof}
First we adopt Holub's proof of Theorem 3.12 in \cite{Ho} to prove that $\{e_i\otimes\cdots\otimes e_i:\, i \in \N\}$ is a $1$-unconditional basic sequence of $\iepx$. Let $I: X \to X$ be the identity operator and for $\theta_i = \pm 1$ ($i \in \N$), define $T: X\to X$ by $T(\si a_ie_i) = \si \theta_ia_ie_i$ for every $x = \si a_ie_i \in X$. Then $\|T\| \lleq 1$. Now for any $m \in \N$, by Lemma 2.1,
\begin{eqnarray*}
\left\|\smi \theta_ia_ie_i\otimes\cdots\otimes e_i\right\|_\epsilon &=& \left\|\smi a_iT(e_i)\otimes I(e_i)\otimes\cdots\otimes I(e_i)\right\|_\epsilon\\
 &\lleq& \left\|\smi a_ie_i\otimes\cdots\otimes e_i\right\|_\epsilon.
\end{eqnarray*}
Thus  $\{e_i\otimes\cdots\otimes e_i:\, i \in \N\}$ is a $1$-unconditional basic sequence of $\iepx$.

Next we show that $Q$ is well-defined and bounded with $\|Q\| \lleq 1$. Take any $u = \sum_{i_1,\dots, i_n}  b_{i_1,\cdots, i_n} e_{i_1}\otimes\cdots\otimes e_{i_n} \in \iepx$. For every $p, q \in \N$ with $p < q$, let 
$$
u_{p,q} = \sum_{i_1,\dots,i_n = p}^q b_{i_1,\cdots, i_n} e_{i_1}\otimes\cdots\otimes e_{i_n}.
$$
Then there exist $x_{j,k} = \si a_{i,j,k}e_i \in X$, $k = 1, \dots, m$ and $j = 1, \dots, n$ such that
$$
u_{p,q} = \smk x_{1,k}\otimes\cdots\otimes x_{n,k}.
$$
Thus
$$
b_{i,\cdots,i} = \smk a_{i,1,k}\cdots a_{i,n,k}, \quad p \lleq i \lleq q.
$$
By Rademacher averaging,
\begin{eqnarray*}
\left\|\sum_{i=p}^q b_{i,\cdots,i} e_i\otimes\cdots\otimes e_i\right\|_{\epsilon} &=& \left\|\smk \sum_{i=p}^q a_{i,1,k}\cdots a_{i,n,k} e_i\otimes\cdots\otimes e_i\right\|_{\epsilon}\\ 
 &=& \left\|\smk\sum_{i=p}^q (a_{i,1,k}e_i)\otimes\cdots\otimes (a_{i,n,k}e_i)\right\|_{\epsilon}\\
 &=& \left\|\smk\int_0^1\Big(\sum_{i=p}^q  a_{i,1,k}r_i(t)e_i\Big)\otimes\cdots\otimes \Big(\sum_{i=p}^q a_{i,n,k}r_i(t)e_i\Big)\,dt\right\|_{\epsilon}\\
 &\lleq& \int_0^1\left\|\smk\Big(\sum_{i=p}^q  a_{i,1,k}r_i(t)e_i\Big)\otimes\cdots\otimes \Big(\sum_{i=p}^q a_{i,n,k}r_i(t)e_i\Big)\right\|_{\epsilon}\,dt\\
 &=& \int_0^1\left\|\smk\Big(T_t(x_{1,k})\Big)\otimes\cdots\otimes \Big(T_t(x_{n,k})\Big)\right\|_{\epsilon}\,dt\\
 &\lleq& \int_0^1\|T_t\|^n\cdot \left\|\smk x_{1,k}\otimes\cdots\otimes x_{n,k}\right\|_{\epsilon}\,dt\\
 &\lleq& \left\|\smk x_{1,k}\otimes\cdots\otimes x_{n,k}\right\|_{\epsilon} = \Big\|u_{p,q}\Big\|_\epsilon,
\end{eqnarray*}
where $T_t: X\to X$ is defined by $T_t(x) = \si a_ir_i(t)e_i$ for every $x = \si a_ie_i \in X$ and every $t \in [0,1]$. Therefore, for every $p,q \in \N$ with $p < q$,
$$
\left\|\sum_{i=p}^q b_{i,\cdots, i} e_i\otimes\cdots\otimes e_i\right\|_{\epsilon} \lleq \left\|\sum_{i_1,\cdots,i_n = p}^q b_{i_1,\cdots, i_n} e_{i_1}\otimes\cdots\otimes e_{i_n}\right\|_{\epsilon},
$$
which implies that $Q$ is well-defined and bounded with $\|Q\| \lleq 1$.
\end{proof}

\vspace{.1in}
\th{Lemma 3.2.}
The tensor diagonal $\{e_i\otimes\cdots\otimes e_i:\, i \in \N\}$ is a $1$-unconditional basis of $\idspx$ and the projection $Q_s: \iespx\rightarrow \idspx$ defined by
$$
Q_s(e_{i_1}\otimes_s\cdots\otimes_s e_{i_n}) = \left\{
\begin{array}{ll}
e_{i_1}\otimes_s\cdots\otimes_s e_{i_n}, & \mbox{if}\;\; i_1 = \cdots = i_n,\\
0, & \mbox{otherwise}.
\end{array}
\right.
$$
is bounded  with $\|Q_s\| \lleq 1$.
\endth

\begin{proof}
For any $m \in \N$, let $u = \smi a_i e_i\otimes\cdots\otimes e_i$, and let $P_u: X^\ast \to \R$ be the $n$-homogeneous polynomial defined in (2.5), that is, 
$$
P_u(x^\ast) = \smi a_i \cdot \Big(x^\ast(e_i)\Big)^n, \quad \forall\; x^\ast \in X^\ast.
\eqno{(3.1)}
$$

For every $a \in \R$ and every $p > 0$ we define $a^p = \mbox{sign}(a) \cdot |a|^p$. Note that $X$ has a 1-unconditional basis $\{e_i: i \in \N\}$. It can be a Banach lattice with the order defined coordinatewise. Also note that $X^\ast$ is a sequence space via $x^\ast \leftrightarrow (x^\ast(e_1), x^\ast(e_2), \dots)$ for every $x^\ast \in E^\ast$. Thus for every $x^\ast, y^\ast \in X^\ast$, from functional calculation, $\big(x^{\ast p} + y^{\ast p}\big)^{\frac{1}{p}}$ is defined (coordinatewise) to be an element of $X^\ast$ (see, e.g., \cite[Section 1.d]{LT}), that is, 
$$
\big(x^{\ast p} + y^{\ast p}\big)^{\frac{1}{p}}(e_i) =  \big(x^\ast(e_i)^p + y^\ast(e_i)^p\big)^{\frac{1}{p}}, \quad i = 1, 2, \dots.
$$
If, moreover, $x^\ast\perp y^\ast$, then  $\big(x^{\ast p} + y^{\ast p}\big)^{\frac{1}{p}} = x^\ast + y^\ast$ by \cite[Lemma 3]{can}. Thus
$$
\big(x^\ast(e_i) + y^\ast(e_i)\big)^{p} =  x^\ast(e_i)^p + y^\ast(e_i)^p, \quad i = 1, 2, \dots.
\eqno{(3.2)}
$$
It follows from (3.1) and (3.2) that $P_u(x^\ast + y^\ast) = P_u(x^\ast) + P_u(y^\ast)$ for every $x^\ast, y^\ast \in X^\ast$ with $x^\ast \perp y^\ast$. Thus $P_u$ is orthogonally additive (see \cite{BB}). Let $T_u: X^\ast\times\cdots\times X^\ast \to \R$ be the $n$-linear operator defined in (2.1), that is,
$$
T_u(x_1^\ast, \dots, x_n^\ast) = \smi a_i \cdot x_1^\ast(e_i)\cdots x_n^\ast(e_i), \quad \forall\; x_1^\ast, \dots, x_n^\ast \in X^\ast.
$$
Then $T_u$ is the symmetric $n$-linear operator associated to $P_u$. It follows from \cite[Theorem 5.4]{BB1} that $\|T_u\| = \|P_u\|$ and hence, $\|u\|_{s,\epsilon} = \|P_u\| = \|T_u\| = \|u\|_{\epsilon}$. 

Now let $\theta_i = \pm 1$ ($i \in \N$). Then by (2.4) and Lemma 3.1,
\begin{eqnarray*}
\left\|\smi \theta_ia_ie_i\otimes\cdots\otimes e_i\right\|_{s,\epsilon} &\lleq& \left\|\smi \theta_ia_ie_i\otimes\cdots\otimes e_i\right\|_{\epsilon}\\
 &\lleq& \left\|\smi a_ie_i\otimes\cdots\otimes e_i\right\|_{\epsilon} = \|u\|_\epsilon\\
 &=& \|u\|_{s,\epsilon} = \left\|\smi a_ie_i\otimes\cdots\otimes e_i\right\|_{s,\epsilon}.
\end{eqnarray*}
Thus  $\{e_i\otimes\cdots\otimes e_i:\, i \in \N\}$ is a 1-unconditional basic sequence of $\iespx$.

Similarly to the proof of Lemma 3.1, we can use the Rademacher averaging formula to show that $Q_s$ is well-defined and bounded with $\|Q_s\| \lleq 1$.
\end{proof}

\vspace{.2in}
By a {\bf Banach lattice with a Schauder basis} we mean a Banach lattice in which the unit vectors form a basis and the order is defined coordinatewise. It follows that such a Schauder basis is 1-unconditional. Conversely, every Banach space with a 1-unconditional basis is a Banach lattice with the order defined coordinatewise. In what follows $E_1, \dots, E_n$ are Banach lattices with (1-unconditional) basis $\{e^{(1)}_i:\, i \in \N\}, \dots, \{e^{(n)}_i:\, i \in \N\}$, respectively.

\vspace{.1in}
\th{Theorem 3.3.}
The tensor basis $\{e^{(1)}_{i_1}\otimes\cdots\otimes e^{(n)}_{i_n}:\, (i_1, \dots, i_n) \in \N^n\}$ with any order is a ($1$-unconditional) basis of $E_1\check{\otimes}_{|\epsilon|} \cdots \check{\otimes}_{|\epsilon|} E_n$.
\endth

\begin{proof}
First, we will show that $(e^{(1)}_{i_1}\otimes\cdots\otimes e^{(n)}_{i_n}) \perp (e^{(1)}_{k_1}\otimes\cdots\otimes e^{(n)}_{k_n})$ provided $(i_1, \dots, i_n) \not= (k_1, \dots, k_n)$. Note that $E_1\check{\otimes}_{|\epsilon|} \cdots \check{\otimes}_{|\epsilon|} E_n$ can be considered as a closed sublattice of $\Lr(E_1^\ast, \dots, E_n^\ast; \R)$. It suffices to show that 
$$
\langle (e^{(1)}_{i_1}\otimes\cdots\otimes e^{(n)}_{i_n}) \wedge (e^{(1)}_{k_1}\otimes\cdots\otimes e^{(n)}_{k_n}), (x_1^\ast, \dots, x_n^\ast) \rangle = 0
$$
for every $x_1^\ast \in E_1^{\ast +}, \dots, x_n^\ast \in E_n^{\ast +}$. Let 
$$
\alpha_1 = x_1^\ast(e^{(1)}_{k_1}), \dots, \alpha_n = x_n^\ast(e^{(n)}_{k_n})\quad \mbox{and}\quad u_{1,1}^\ast = \alpha_1 f^{(1)}_{k_1}, \dots, u_{n,1}^\ast = \alpha_n f^{(n)}_{k_n},
$$
where $\{f^{(1)}_i:\, i \in \N\}, \dots, \{f^{(n)}_i:\, i \in \N\}$ are, respectively, the appropriate bi-orthogonal functionals on $E_1^\ast, \dots, E_n^\ast$. Then for every $x = \si b_i e_i^{(1)} \in E_1^{+}$,
$$
u_{1,1}^\ast(x) = \alpha_1f^{(1)}_{k_1}(x) = x_1^\ast(e^{(1)}_{k_1}) b_{k_1} \lleq \si b_i x_1^\ast(e^{(1)}_{i}) = x_1^\ast(x).
$$
It follows that $u_{1,1}^\ast \lleq x_1^\ast$ and similarly, $u_{2,1}^\ast \lleq x_2^\ast, \dots, u_{n,1}^\ast \lleq x_n^\ast$. Let $u^\ast_{1,2} = x^\ast_1 - u^\ast_{1,1}, \dots, u^\ast_{n,2} = x^\ast_n - u^\ast_{n,1}$. Then $(u^\ast_{1,1}, u^\ast_{1,2}), \dots, (u^\ast_{n,1}, u^\ast_{n,2})$ are partitions of $x^\ast_1, \dots, x^\ast_n$ respectively. By \cite[Proposition 2.2]{BBL},
\begin{eqnarray*}
 &&\langle (e^{(1)}_{i_1}\otimes\cdots\otimes e^{(n)}_{i_n}) \wedge (e^{(1)}_{k_1}\otimes\cdots\otimes e^{(n)}_{k_n}), (x_1^\ast, \dots, x_n^\ast) \rangle\\
 &\lleq& \sum^2_{j_1,\dots,j_n = 1}\langle e^{(1)}_{i_1}\otimes\cdots\otimes e^{(n)}_{i_n}, (u_{1,j_1}^\ast, \dots, u_{n,j_n}^\ast) \rangle\,\wedge\,\langle e^{(1)}_{k_1}\otimes\cdots\otimes e^{(n)}_{k_n}, (u_{1,j_1}^\ast, \dots, u_{n,j_n}^\ast) \rangle.
\end{eqnarray*}
If all $j_m$'s are 1, then the general term 
\begin{eqnarray*}
 && \langle e^{(1)}_{i_1}\otimes\cdots\otimes e^{(n)}_{i_n}, (u_{1,j_1}^\ast, \dots, u_{n,j_n}^\ast) \rangle\,\wedge\,\langle e^{(1)}_{k_1}\otimes\cdots\otimes e^{(n)}_{k_n}, (u_{1,j_1}^\ast, \dots, u_{n,j_n}^\ast) \rangle\\
 &\lleq& \langle e^{(1)}_{i_1}\otimes\cdots\otimes e^{(n)}_{i_n}, (u_{1,1}^\ast, \dots, u_{n,1}^\ast) \rangle
= u_{1,1}^\ast(e^{(1)}_{i_1}) \cdots u_{n,1}^\ast(e^{(n)}_{i_n})\\
 &=& \alpha_1 f^{(1)}_{k_1}(e^{(1)}_{i_1}) \cdots \alpha_n f^{(n)}_{k_n}(e^{(n)}_{i_n}) = 0
\end{eqnarray*}
since $(i_1, \dots, i_n) \not= (k_1, \dots, k_n)$. If at least one of $j_m$'s is 2, say $j_1 = 2$, then the general term 
\begin{eqnarray*}
 && \langle e^{(1)}_{i_1}\otimes\cdots\otimes e^{(n)}_{i_n}, (u_{1,j_1}^\ast, \dots, u_{n,j_n}^\ast) \rangle\,\wedge\,\langle e^{(1)}_{k_1}\otimes\cdots\otimes e^{(n)}_{k_n}, (u_{1,j_1}^\ast, \dots, u_{n,j_n}^\ast) \rangle\\
 &\lleq& \langle e^{(1)}_{k_1}\otimes e^{(2)}_{k_2}\otimes\cdots\otimes e^{(n)}_{k_n}, (u_{1,2}^\ast, u^\ast_{2,j_2}, \dots, u_{n,j_n}^\ast) \rangle\\
 &=& u_{1,2}^\ast(e^{(1)}_{k_1})\cdot u_{2,j_2}^\ast(e^{(2)}_{k_2}) \cdots u_{n,j_n}^\ast(e^{(n)}_{k_n})\\
 &=& \Big(\alpha_1 - \alpha_1 f^{(1)}_{k_1}(e^{(1)}_{k_1}) \Big) \cdot u_{2,j_2}^\ast(e^{(2)}_{k_2}) \cdots u_{n,j_n}^\ast(e^{(n)}_{k_n}) = 0.
\end{eqnarray*}
Therefore, 
$$
\langle (e^{(1)}_{i_1}\otimes\cdots\otimes e^{(n)}_{i_n}) \wedge (e^{(1)}_{k_1}\otimes\cdots\otimes e^{(n)}_{k_n}), (x_1^\ast, \dots, x_n^\ast) \rangle = 0
$$ 
and hence, $(e^{(1)}_{i_1}\otimes\cdots\otimes e^{(n)}_{i_n}) \perp (e^{(1)}_{k_1}\otimes\cdots\otimes e^{(n)}_{k_n})$. 

Being a disjoint sequence in a Banach lattice, $\{e^{(1)}_{i_1}\otimes\cdots\otimes e^{(n)}_{i_n}:\, (i_1, \dots, i_n) \in \N^n\}$ is a $1$-unconditional basic sequence of $E_1\check{\otimes}_{|\epsilon|} \cdots \check{\otimes}_{|\epsilon|} E_n$. It is left to show that its span is dense in $E_1\check{\otimes}_{|\epsilon|} \cdots \check{\otimes}_{|\epsilon|} E_n$.

Take any $x_i \in E_i$ with $\|x_i\| \lleq 1$ for $i = 1, \dots, n$. Given any $\sigma \in (0,1)$, we can find basis projections $P_i$ on $E_i$, respectively, such that 
$$
y_i = P_i(x_i) \quad\mbox{and}\quad \|y_i - x_i\| \lleq \sigma, \quad i = 1, \dots, n.
$$
Then
\begin{eqnarray*}
\|x_1\otimes x_2 - y_1\otimes y_2\|_{|\epsilon|} &=& \|x_1\otimes (x_2 - y_2) + (x_1 - y_1)\otimes y_2\|_{|\epsilon|}\\
 &\lleq& \|x_1\| \cdot \|x_2 - y_2\| + \|x_1 - y_1\| \cdot \|y_2\| \lleq 2\sigma,
\end{eqnarray*}
and similarly,
$$
\|x_1\otimes\cdots\otimes x_n - y_1\otimes\cdots\otimes y_n\|_{|\epsilon|} \lleq n\sigma.
$$
Since $y_1\otimes\cdots\otimes y_n$ is in the span $\{e^{(1)}_{i_1}\otimes\cdots\otimes e^{(n)}_{i_n}:\, (i_1, \dots, i_n) \in \N^n\}$, it follows that $x_1\otimes\cdots\otimes x_n$ is in  the closure of this span. Thus this span is dense in $E_1{\otimes}\cdots{\otimes}E_n$ and hence, dense in $E_1\check{\otimes}_{|\epsilon|} \cdots \check{\otimes}_{|\epsilon|} E_n$.
\end{proof}

\vspace{.1in}
In what follows $E$ is a Banach lattice with a (1-unconditional) basis $\{e_i:\, i \in \N\}$. The following consequence comes straightforward from the previous theorem.

\vspace{.1in}
\th{Corollary 3.4.}
The tensor basis $\{e_{i_1}\otimes\cdots\otimes e_{i_n}:\, (i_1, \dots, i_n) \in \N^n\}$ with any order is a ($1$-unconditional) basis of $\iep$, and the tensor basis $\{e_{i_1}\otimes_s\cdots\otimes_s e_{i_n}:\, (i_1, \dots, i_n) \in \N^n, i_1 \ggeq\cdots\ggeq i_n\}$ with any order is a ($1$-unconditional) basis of $\iesp$.
\endth

\vspace{.1in}
Let $\idpp$ (resp. $\idspp$) denote the {\it main diagonal space} of $\iep$ (resp. $\iesp$), that is, the closed subspace spanned in $\iep$ (resp. in $\iesp$) by the {\it tensor diagonal} $\{e_i\otimes\cdots\otimes e_i:\, i \in \N\}$. It follows from the above lemma that $\{e_i\otimes\cdots\otimes e_i:\, i \in \N\}$ is a (1-unconditional) basis of both $\idpp$ and $\idspp$. 

Recall that we already have other two main diagonal spaces $\idpe$ and $\idsp$ introduced at the beginning of this section. Next we will show that all four main diagonal spaces are pairwise isometrically isomorphic. First we need the following H\"{o}lder inequality and then a lemma.

\vspace{.1in}
\th{H\"{o}lder Inequality.}
Let $1 \lleq p_1, \dots, p_n \lleq \infty$ such that $1/p_1 + \cdots + 1/p_n = 1$. Then for every $(b^{(k)}_i)_{i} \in \ell_{p_k}$, $k = 1, \dots, n$, we have 
$$
\si \left|b_i^{(1)}\cdots b_i^{(n)}\right| \lleq \left\|\big(b^{(1)}_i\big)_i\right\|_{\ell_{p_1}} \cdots \left\|\big(b^{(n)}_i\big)_i\right\|_{\ell_{p_n}}.
$$
\endth

\vspace{.1in}
\th{Lemma 3.5.}
For $x^\ast_1, \dots, x^\ast_n \in E^{\ast+}$, define $x^\ast: E\to \R$ by
$$
x^\ast(x) = \si a_i\big(x^\ast_1(e_i)\big)^{\frac{1}{n}} \cdots \big(x^\ast_n(e_i)\big)^{\frac{1}{n}}
$$
for every $x = \si a_ie_i \in E$. Then $x^\ast \in E^{\ast+}$ with $\|x^\ast\| \lleq \big\|x^\ast_1\big\|^{\frac{1}{n}} \cdots \big\|x^\ast_1\big\|^{\frac{1}{n}}$.
\endth

\begin{proof}
For any $m \in \N$, by H\"{o}lder inequality,
\begin{eqnarray*}
 & & \smi \Big|a_i\big(x^\ast_1(e_i)\big)^{\frac{1}{n}} \cdots \big(x^\ast_n(e_i)\big)^{\frac{1}{n}}\Big|\\
 &=& \smi \Big|a_ix^\ast_1(e_i)\Big|^{\frac{1}{n}} \cdots \Big|a_ix^\ast_n(e_i)\Big|^{\frac{1}{n}}\\
 &\lleq& \left(\smi|a_ix^\ast_1(e_i)|\right)^{\frac{1}{n}} \cdots \left(\smi|a_ix^\ast_n(e_i)|\right)^{\frac{1}{n}}\\
 &=& \left(x^\ast_1\big(\smi \pm a_ie_i\big)\right)^{\frac{1}{n}} \cdots \left(x^\ast_n\big(\smi \pm a_ie_i\big)\right)^{\frac{1}{n}}\\ 
 &\lleq& \big\|x^\ast_1\big\|^{\frac{1}{n}} \cdots \big\|x^\ast_n\big\|^{\frac{1}{n}} \cdot \left\|\smi a_ie_i\right\|.
\end{eqnarray*}
Thus $x^\ast$ is well defined and $x^\ast \in E^{\ast+}$ with $\|x^\ast\| \lleq \big\|x^\ast_1\big\|^{\frac{1}{n}} \cdots \big\|x^\ast_1\big\|^{\frac{1}{n}}$.
\end{proof}

\vspace{.1in}
\th{Theorem 3.6.}
All four main diagonal spaces $\idpe$, $\idsp$, $\idpp$, and $\idspp$ are pairwise isometrically isomorphic.
\endth

\begin{proof}
First we show that $\idpp$ is isometrically isomorphic to $\idspp$. Since $\{e_i\otimes\cdots\otimes e_i: \; i \in \N\}$ is a basis of both $\idpp$ and $\idspp$, it suffices to show that $\|u\|_{|\epsilon|} = \|u\|_{s,|\epsilon|}$ for every $u = \smi a_i e_i\otimes\cdots\otimes e_i \in \e$. Without loss of generality, we assume that $u \ggeq 0$, that is, $a_i \ggeq 0$ for $ i = 1, \dots, m$. For any $\sigma > 0$ there exist $x^\ast_1, \dots, x^\ast_n \in B_{E^{\ast+}}$ such that
$$
\|u\|_{|\epsilon|} \lleq \smi a_i x^\ast_1(e_i) \cdots x^\ast_n(e_i) + \sigma.
$$
Let $x^\ast$ be defined in Lemma 3.5. Then $x^\ast \in B_{E^{\ast+}}$ and for each $i \in \N$,
$$
x^\ast(e_i) = \big(x^\ast_1(e_i)\big)^{\frac{1}{n}} \cdots \big(x^\ast_n(e_i)\big)^{\frac{1}{n}}.
$$
Thus 
$$
\|u\|_{s,|\epsilon|} \ggeq \smi a_i\big(x^\ast(e_i)\big)^{n} = \smi a_i x^\ast_1(e_i) \cdots x^\ast_n(e_i) \ggeq \|u\|_{|\epsilon|} - \sigma.
$$
It follows that $\|u\|_{s,|\epsilon|} \ggeq \|u\|_{|\epsilon|}$. From their definitions it is trivial that $\|u\|_{s,|\epsilon|} \lleq \|u\|_{|\epsilon|}$ and hence, $\|u\|_{s,|\epsilon|} = \|u\|_{|\epsilon|}$.

Secondly, we show that $\idpe$ is isometrically isomorphic to $\idpp$. Since $\{e_i\otimes\cdots\otimes e_i: \; i \in \N\}$ is a basis of both $\idpe$ and $\idpp$, it suffices to show that $\|u\|_{\epsilon} = \|u\|_{|\epsilon|}$ for every $u = \smi a_i e_i\otimes\cdots\otimes e_i \in \e$. By Lemma 3.1,
\begin{eqnarray*}
\big\|u\big\|_{|\epsilon|} &=& \big\||u|\big\|_{|\epsilon|} = \left\|\smi |a_i| e_i\otimes\cdots\otimes e_i\right\|_{|\epsilon|}\\
 &=& \left\|\smi |a_i| e_i\otimes\cdots\otimes e_i\right\|_{\epsilon} = \left\|\smi \pm a_i e_i\otimes\cdots\otimes e_i\right\|_{\epsilon}\\
 &\lleq& \left\|\smi a_i e_i\otimes\cdots\otimes e_i\right\|_{\epsilon} = \big\|u\big\|_{\epsilon}.
\end{eqnarray*}
On the other hand, it follows from their definitions that $\|u\|_\epsilon \lleq \|u\|_{|\epsilon|}$ and hence, $\|u\|_\epsilon = \|u\|_{|\epsilon|}$.

Finally we show that $\idpe$ is isometrically isomorphic to $\idsp$. Since $\{e_i\otimes\cdots\otimes e_i: \; i \in \N\}$ is a basis of both $\idpe$ and $\idsp$, it suffices to show that $\|u\|_{\epsilon} = \|u\|_{s,\epsilon}$ for every $u = \smi a_i e_i\otimes\cdots\otimes e_i \in \e$, which was proved in the proof of Lemma 3.2.
\end{proof}

\vspace{.2in}
\bibliographystyle{amsplain}

\begin{thebibliography}{99}


\bibitem{BB}
Q. Bu and G. Buskes, {Polynomials on Banach lattices and positive tensor products}, {\it J. Math. Anal. Appl.} {\bf 388} (2012), 845--862.

\bibitem{BB1}
Q. Bu and G. Buskes, Diagonals of projective tensor products and orthogonally additive polynomials, {\it Studia Math.} {\bf 201} (2014), 1--15.

\bibitem{BBL}
Q. Bu, G. Buskes, and Y. Li, AM- and AL-spaces of polynomials on Banach lattices, {\it Proc. Edinburgh Math. Soc.},  

\bibitem{can}
Q. Bu, G. Buskes, A.I. Popov, A. Tcaciuc, and V.G. Troitsky, {The 2-concavification of a Banach lattice equals the diagonal of the Fremlin tensor square}, {\it Positivity} {\bf 17} (2012), 283--298.

\bibitem{Ca}
F.Cabello-Sanchez, {Complemented subspaces of spaces of multilinear forms and tensor products}, {\it J. Math. Anal. Appl.} {\bf 254} (2001), 645--653.


\bibitem{CG}
D. Carando and D. Galicer, Unconditionality in tensor products and ideals of polynomials, multilinear forms and operators, {\it Quart. J. Math.} {\bf 62} (2011), 845--869.

\bibitem{DF}
A. Defant and K. Floret, Tensor Norms and Operator Ideals, North-Holland, Amesterdam, 1993.

\bibitem{DDGM}
A. Defant, J.C. D\'{i}az, D. Garcia, and M. Maestre, Unconditional basis and Gordon-Lewis constants for spaces of polynomials, {\it J. Funct. Anal.} {\bf 182} (2001), 119--145.

\bibitem{Di}
S. Dineen, {Complex Analysis on Infinite Dimensional Spaces}, Springer, 1999.

\bibitem{Fl}
K. Floret, {Natural norms on symmetric tensor products of normed spaces}, {\it Note Mat.} {\bf 17} (1997), 153--188.


\bibitem{Ge}
B.R. Gelbaum and  J.Gil de Lamadrid, {Bases of tensor products of Banach spaces}, {\it Pacific J. Math.} {\bf 11} (1961), 1281-1286.


\bibitem{GR1}
B.C. Grecu and R.A. Ryan, {Schauder bases for symmetric tensor products}, {\it Publ. RIMS, Kyoto Univ.} {\bf 41} (2005), 459--469. 

\bibitem{Ho}
J.R. Holub, {Tensor product bases and tensor diagonals}, {\it Trans. Amer. Math. Soc.} {\bf 151} (1970), 563-579.


\bibitem{KP}
S. Kwapie\'{n} and A. Pelczy\'{n}ski, {The main triangle projection in matrix spaces and its applications}, {\it Studia Math.}  {\bf 34} (1970), 43--68.

\bibitem{LT}
J. Lindenstrauss and L. Tzafriri, {Classical Banach Spaces II-Function Spaces}, Springer, 1979.


\bibitem{Mu}
J. Mujica, Complex Analysis in Banach Spaces, North-Holland-Math. Stud., {\bf 120}, 1986.

\bibitem{PV}
D. P\'{e}rez-Garc\'{i}a and I. Villanueva, Unconditional bases in tensor products of Hilbert spaces, {\it Math. Scand.} {\bf 96} (2005), 280--288.


\bibitem{Pi}
G. Pisier, Some result on Banach spaces without unconditional structure, {\it Compositio Math.} {\bf 37} (1978), 3--9.

 
\bibitem{R}
R.A. Ryan, {Introduction to Tensor Products of Banach Spaces}, Springer, 2002.

 
\bibitem{Sch}
C. Sch\"{u}tt, Unconditionality in tensor products, {\it Israel J. Math.} {\bf 31} (1978), 209--216.

\end{thebibliography}

\end{document}